\documentclass[format=acmsmall, screen=true, review=false, natbib=true,
anonymous=false, authorversion=false, timestamp=false, authordraft=false]{acmart}

\usepackage{graphicx}
\usepackage{amsmath}
\usepackage{subcaption}
\usepackage[inline]{enumitem}
\usepackage[normalem]{ulem}
\usepackage{booktabs}
\usepackage{array}
\usepackage{multirow}
\begin{teaserfigure}
	\includegraphics[width=\textwidth]{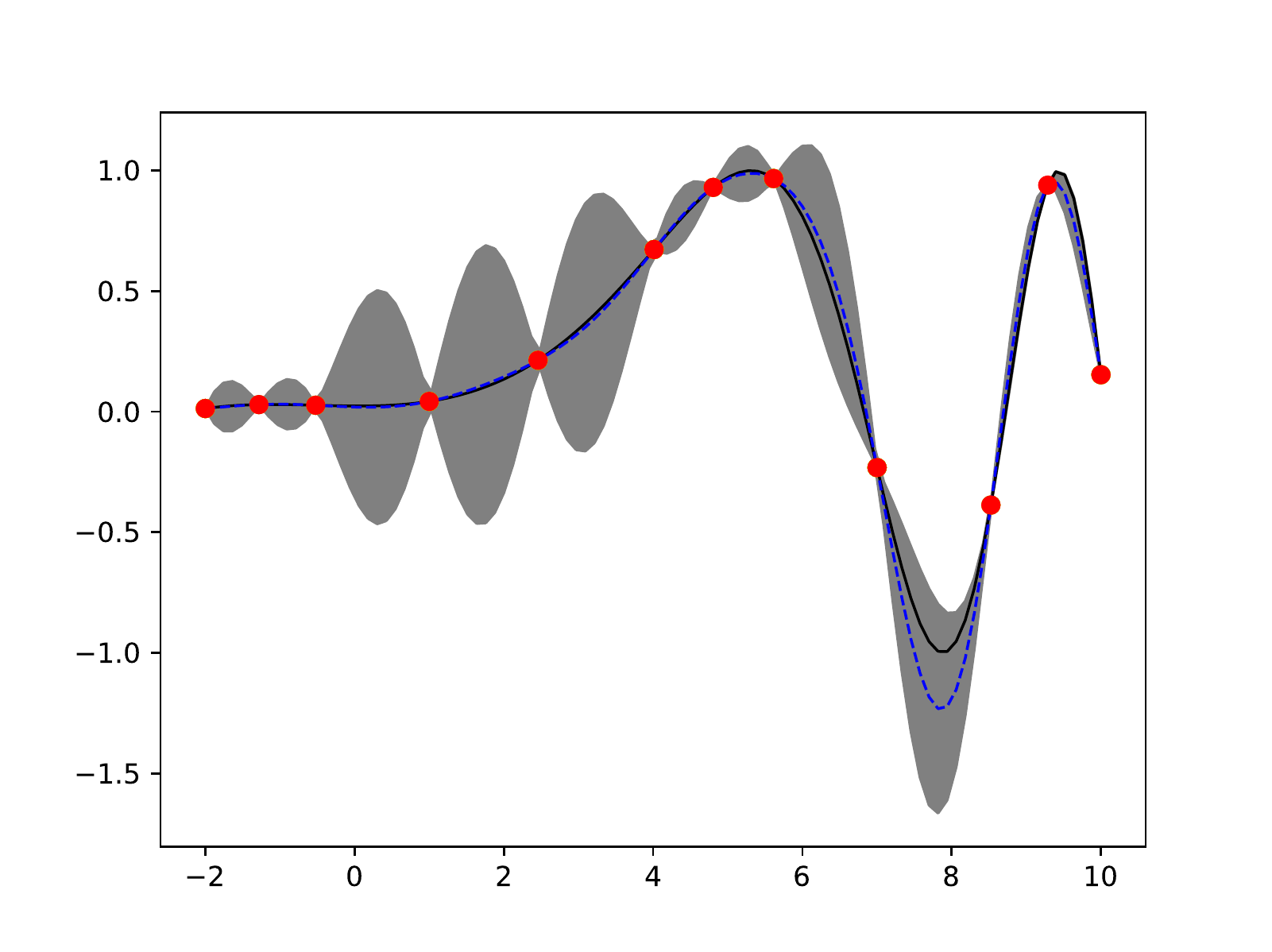}
\caption{Bayesian optimization and Kriging are both enabled by the presented software}
\label{teaser}
\end{teaserfigure}

\acmVolume{1}
\acmNumber{1}
\acmArticle{1}
\acmYear{2018}
\acmMonth{1}
\copyrightyear{1}
\acmJournal{TOMS}

\acmDOI{0000001.0000001}

\received{February 2018}
\received[revised]{1}
\received[accepted]{1}

\usepackage[ruled]{algorithm2e} 

\SetAlFnt{\small}
\SetAlCapFnt{\small}
\SetAlCapNameFnt{\small}
\SetAlCapHSkip{0pt}
\IncMargin{-\parindent}

\newcommand{\mycomment}[1]{} 

\newcommand{\mydraft}[1]{} 

\setcitestyle{numbers,sort&compress}
\begin{document}

\title[PARyOpt]{PARyOpt: A software for  \underline{P}arallel \underline{A}synchronous \underline{R}emote Ba\underline{y}esian \underline{Opt}imization}

\author{Balaji Sesha Sarath Pokuri}
\orcid{0000-0002-5816-0184}
\affiliation{Department of Mechanical Engineering, Iowa State University}
\email{balajip@iastate.edu}

\author{Alec Lofquist}
\orcid{}
\affiliation{Department of Computer Engineering, Iowa State University}
\email{lofquist@iastate.edu}

\author{Chad Risko}
\orcid{}
\affiliation{Department of Chemistry, University of Kentucky}
\email{chad.risko@uky.edu}

\author{Baskar Ganapathysubramanian}
\orcid{}
\affiliation{Department of Mechanical Engineering, Iowa State University}
\email{baskarg@iastate.edu}
\authornote{Corresponding author}

\begin{abstract}
	PARyOpt
	 is a python based implementation of the Bayesian optimization
	routine designed for remote and asynchronous function evaluations. Bayesian
	optimization is especially attractive for computational optimization due to
	its \emph{low cost function footprint} as well as the ability to account
	for uncertainties in data. A key challenge to efficiently deploy any optimization strategy on distributed computing systems is the synchronization step, where data from multiple function calls is assimilated to identify the next campaign of function calls. Bayesian optimization provides an elegant approach to overcome this issue via asynchronous updates. We formulate, develop and implement a parallel, asynchronous variant of Bayesian optimization. The framework is robust and resilient to external failures. We show how such asynchronous evaluations help reduce the total optimization wall clock time for a suite of test problems. Additionally, we show how the software design of the framework allows easy extension to response surface reconstruction (Kriging), providing a high performance software for autonomous exploration. The software is available on PyPI, with examples and documentation. 

\end{abstract}

\maketitle

\section{Introduction}

{\it Simulation-based design} is a resource efficient approach of solving inverse and/or design problems in science and engineering. The goal of inverse/design problems is to identify conditions -- these can be initial conditions, boundary conditions or property/coefficient distribution -- that result in a desired behavior of the engineered system. Examples of these are ubiquitous. A good example is in the identification of tailored processing conditions that result in electronic devices with high performance metrics. In most electronics manufacturing, it has been shown that varying processing conditions can critically impact device properties and the identification of optimal processing conditions is a key problem from the financial and sustainability standpoint\cite{zhao2016vertical, pokuri2017nanoscale, wodo2012modeling}. Another example is in the identification/design of useful molecular architectures that exhibit a suite of desired physicochemical properties (absorption, miscibility, solubility, toxicity, among others). This has substantial implications in the pharmaceutical and the chemical industries. In all such cases of simulation-based engineering, considerable effort has been expended to construct excellent "forward" models of the engineering problem, i.e. models that map the set of input conditions, boundary conditions and property distributions to the output quantity of interest ($\mathcal{F}:\texttt{input}\rightarrow \texttt{property}$). The design problem, however, calls for knowledge of the {\bf reverse mapping} ($\mathcal{F}^{-1}:\texttt{property}\rightarrow \texttt{input}$ ) that maps a desired value of the output to a set of inputs. 

Explicit construction of the reverse mapping is unfeasible in most applications with complex forward models. The design problem of identification is usually posed as an optimization problem. That is, the forward model, $\mathcal{F}$, is solved multiple times within an optimization framework to identify those input conditions that minimize a cost function, with the (argument of the) minima representing the desired input values. In complex engineering applications, this cost function is often very tedious to calculate. It can involve simulations which are computationally expensive. In such situations, these `forward' simulations are performed on large high performance computing (HPC) resources. This federated approach to simulations is quite common, and
gives several benefits such as parallelization and cost effectiveness. However, such large shared compute resources often work with a job scheduler that balances the load across several processes, users and projects, leading to uncertain times of completion (queueing time + initialization time + compute time) for each simulation. Added to that is the logistical effort of managing
data between a local machine and a HPC client. Therefore, in such applications, one can expect the cost function to:
\begin{enumerate*}[label=\itshape\alph*\upshape)]
	\item be solved on a HPC cluster,
	\item have long, uncertain compute times,
	\item require some form of data transfer across computers, and
	\item fail with a small (but finite) probability due to hardware/network
		failures.
\end{enumerate*}
In recent years, there have been various approaches proposed to solve optimization problems under such constraints. 

Bayesian Optimization (BO) is a one such approach, especially for expensive cost function
optimization. BO works by adaptively sampling the input parameters to construct a surrogate of the forward model through basis expansion. By efficiently constructing a surrogate, it achieves
two purposes: one, to interpolate the surrogate through all existing data
points and two, to give confidence estimates on the constructed surrogate. Furthermore, this
information provides a rigorous way to determine the next set of evaluations/
experiments (the adaptive sampling). This form of surrogate construction helps to directly reduce the total number of function evaluations and iterations, helping to efficiently
identify the optima. Additionally, the Gaussian (mean, variance) description used for construction of the surrogate enables a diversity of searches, ranging from pure
exploration to pure exploitation. Exploration helps to build confidence in the
constructed surrogate and exploitation takes advantage of this confidence to
find optima. 
A mixture of exploration and exploitation avoids the optimizer to be stuck at a local optima. Hence, efficient BO usage involves strategic, possibly simultaneous, exploration and
exploitation searches. This parallelized cost function evaluation can often increase the total time of optimization, often due to delays in resource availability and management. BO is superior to other algorithms in this aspect, since the update of the surrogate is Gaussian and hence the order of data assimilation does not affect the posterior surrogate. Additionally, asynchronous data assimilation~\cite{Janusevskis2012} can also be incorporated into this algorithm, with unusually long function evaluations (could be due to difficult convergence, hardware failures and several other logistic reasons) being assimilated in later iterations. Another common challenge with
integration with an HPC system is to harmonize the several runs across simulation nodes. One may also be interested to post-process simulation results to calculate the actual cost function value. That will require a complete setup  to manage data across several systems, quite often using a secure shell protocol.

\mycomment{<how do we give these advantages through PARyOpt>}

Utilizing these advantages and addressing the concerns, we present a
parallelized implementation of Bayesian Optimization that can:
\begin{enumerate*}[label=\itshape\alph*\upshape)]
	\item perform asynchronous updates per iteration,
	\item comes with a secure shell login (SSH) module for HPC integration, and
	\item has fault tolerant restart capability.
\end{enumerate*}
Using PARyOpt, the user can thus perform optimization on a local machine, with
the flexibility to perform simulations either on the same local machine or a
remote machine with/without a job scheduler. To tackle delays in simulation
time, we provide an asynchronous evaluator class that can manage and keep track
of status of simulations and adaptively assimilate them with completion. To
tackle hardware/software failures, we provide functionality to restart
optimization in case of head-node failure and also functionality to re-evaluate
cost function in case of simulation failures. In the presented test cases, the
asynchronous remote evaluator has improved the total time of optimization by up
to $50\%$. Additionally, the software can also be used for general surface
reconstruction, commonly known as Kriging.

The rest of the paper is structured as follows -- in section
\ref{bayes_discussion}, we discuss the basic math behind surrogate
construction, Bayesian update and acquisition function, and potential avenues of
parallelization. Section~\ref{async} discusses current implementation of the
algorithm, and development and usage of the asynchronous evaluator class. In
section \ref{results}, several standard optimization problems as well as
kriging are considered along with a thorough analysis of the novel asynchronous
evaluator class. Finally in section~\ref{future}, code availability,
reproducibility and future directions of this software are discussed.

\section{Bayesian update -- Gaussian processes}
\label{bayes_discussion}

In this section, we present a brief discussion of the algorithm of Bayesian
Optimization. The core algorithm is well studied and not the topic of this
work. We refer the reader to \cite{Brochu2010, shahriari2016taking} for a 
detailed explanation of the algorithm, its requirements and limitations.

We consider a general minimization problem:
\begin{equation}
	\min_\mathbf{x} \, y(\mathbf{x})
\end{equation}

Bayesian optimization proceeds through construction of a surrogate cost
function $\tilde{y}(\mathbf{x})$. This surrogate is represented via a basis
function expansion, around each evaluated point~( $\mathbf{x_i}, i=1,... ,N$ ).
This ensures that the surrogate passes through (interpolates) the evaluated
points. In the case of evaluations with noisy data, the surrogate shall pass
within \textit{one} standard deviation from the mean at the evaluated points.
Analytically, the surrogate $\tilde{y}(\mathbf{x})$ after $N$ function
evaluations is represented as

\begin{equation}
	\tilde{y}(\mathbf{x}) = \sum_{i=1,2,..N}w_i\,k(\mathbf{x}_i,\mathbf{x})
\end{equation}

Typically, $k(\mathbf{x}_i,\mathbf{x})$ is a kernel function
(section~\ref{kernels}), i.e., it takes in two arguments, $\mathbf{x}, \,
\mathbf{x}_i$, and returns a scalar representative of correlation of the
function $y(\mathbf{x})$ between the points $\mathbf{x}$ and $\mathbf{x}_i$. The weights $w_i$
are calculated by solving the system of $N$ linear equations in $w_i$. In
matrix notation, this is represented using a covariance matrix~($\mathbf{K}$):

\begin{align}
	\mathbf{K}\, \bar{w} & = y \\
	\mathbf{K}_{i,j} & = k(\mathbf{x}_i,\mathbf{x}_j) , \, \, i,j\in[1,N] \\
	y_i & = y(\mathbf{x}_i) , \, \,i\in[1,N] \\
	\bar{w} & = \{w_i\}, \, \,i\in[1,N]
\end{align}

Hence the weights are calculated through the inversion $\bar{w} =
\mathbf{K}^{-1}\,y$. Note that the covariance matrix $\mathbf{K}$ is a Gram
matrix of a positive definite kernel function, making it symmetric and positive
semi-definite (see section~\ref{kernels}). Furthermore, since with every
iteration only a finite number of rows are added to the covariance matrix,
efficient inversion is possible through incremental Cholesky
decomposition~\cite{Polok2013}. The mean and variance of the surrogate are then
calculated as:

\begin{align}
	\mu(\mathbf{x}_{N+1}) & = \mathbf{k}^T \mathbf{K}^{-1} y_{1:N} \\
	\sigma^2(\mathbf{x}_{N+1}) & = k(\mathbf{x}_{N+1},\mathbf{x}_{N+1}) - \mathbf{k}^T\,\mathbf{K}^{-1}\,\mathbf{k}
\end{align}

where

\begin{equation}
	\mathbf{k} = \mathbf{k}(\mathbf{x}_{N+1}) = [k(\mathbf{x}_1,\mathbf{x}_{N+1})\, k(\mathbf{x}_2,\mathbf{x}_{N+1})\, . . . 	k(\mathbf{x}_N,\mathbf{x}_{N+1})]
\end{equation}

At each iteration, the surrogate is updated with new data from the cost
function. The locations where the next evaluation is done is determined through
optimization of an \emph{acquisition function} (section~\ref{acqs}). An
acquisition function is a means to estimate the new information content at a
location. It uses the mean and variance calculated in the above steps. In total,
the algorithm is summarized in the following figure (figure~\ref{bayes_opt_algo}).

\begin{figure}
	\includegraphics[width=\textwidth]{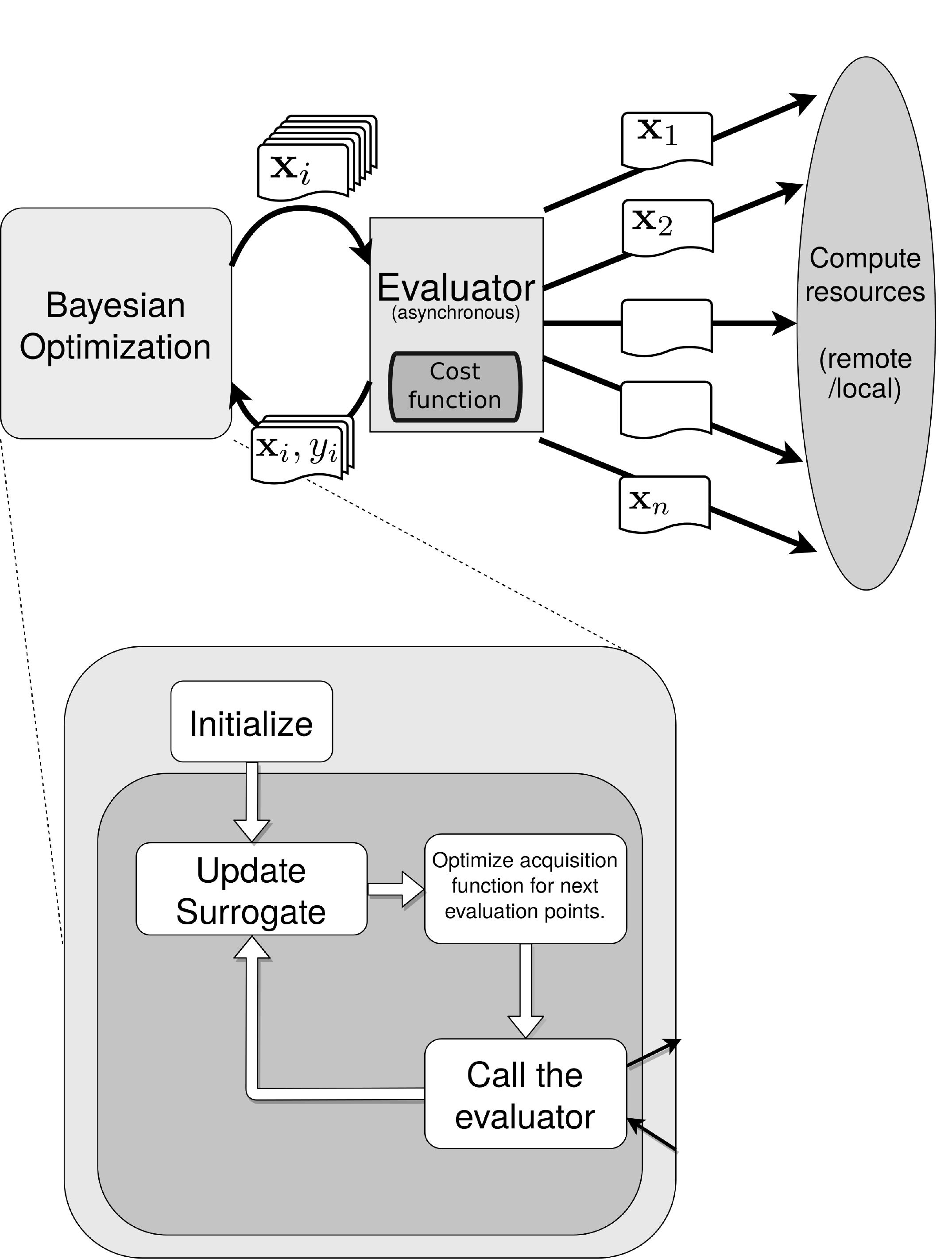}
	\caption{Modular design of bayesian optimization. Note how the evaluator
		can be designed with custom capabilities, such as asynchronocity and remote
		evaluation. The next figure will explain the implementation of an asynchronous
	evaluator.}
	\label{bayes_opt_algo}
\end{figure}

\subsection{Kernel function}
\label{kernels}

\mydraft{
	\begin{itemize}
		\item requirements --
		\item commonly used kernel functions, options available in current
			implementation
	\end{itemize}
}

The kernel function, $f_i(\mathbf{x})=k(\mathbf{x},\mathbf{x}_i)$, is the most
important component of Bayesian Optimization. It embeds the information about
continuity, periodicity and correlation lengths of the underlying cost
function, $y$. In order to be used as a kernel function in Gaussian processes the
function should be a positive semi-definite function, i.e,

\begin{equation}
	\int k(\mathbf{x}, \mathbf{x'})g(\mathbf{x})g(\mathbf{x'})\,d\mu(\mathbf{x})\,d\mu(\mathbf{x'}) \geq 0
\end{equation}

for all $g \in L(\chi, \mu)$, $\mu$ is any measure. An exhaustive description of
the requirements of covariance functions can be found in~\cite{CarlEdward2004}.
Some examples of standard stationary covariance functions include:
\begin{enumerate}
	\item Squared Exponential covariance function
		\begin{equation}
			k_{SE}(r) = \exp(-\frac{r^2}{l^2})
		\end{equation}
	\item Matern class of covariance functions
		\begin{equation}
			k_{Matern}(r) = \frac{2^{1-\nu}}{\Gamma(\nu)}(\frac{\sqrt{2\nu}r}{l})^\nu K_\nu(\frac{\sqrt{2\nu}r}{l})
		\end{equation}
		where $K_\nu$ is the modified Bessel function, $\nu,l$ are positive constants
	\item Exponential covariance function
		\begin{equation}
			k_{E}(r) = \exp(-\frac{r}{l})
		\end{equation}
	\item $\gamma$-exponential function
		\begin{equation}
			k_{gE}(r) = \exp(-(\frac{r}{l})^\gamma)
		\end{equation}
	\item Rational quadratic covariance function
		\begin{equation}
			k_{RQ}(r) = (1 + \frac{r^2}{2\alpha l^2})^{-\alpha}
		\end{equation}
	\item Peicewise polynomials with compact support: Several algorithms are
		discussed in~\cite{wendland2005scattered}, one such example with compact support
		in $D$ dimensions is given by
		\begin{equation}
			k_{ppD,0} = (1-r)^j_+
		\end{equation}where, $j = \lfloor\frac{D}{2}\rfloor + q + 1 $
\end{enumerate}
where $$ r = ||\mathbf{x}-\mathbf{x}'|| $$
While several classes of kernel functions can be used, in our implementation,
we focus on stationary functions (radial basis functions). It contains some of
the most commonly used kernel functions like squared exponential, Matern $3/2$
and Matern $5/2$. It is implemented as a separate kernel function class, so
that users can derive this class and supply their own custom kernel functions.
This class also contains the functionality to set the length scale parameters.
This functionality is used to estimate and tune the current length scales for
the construction of a better surrogate with a given class of kernel functions.

\subsubsection{Maximum Likelihood estimate}

The Maximum Likelihood Estimate (MLE) for a Gaussian process gives information
about the goodness of the surrogate based on the current available data. This
provides a rigorous means of estimating the length scale in the kernel
function. Generally known as hyper-parameter estimation, this is done by
minimizing the MLE of the surrogate~\cite{mongillo2011choosing}. The MLE is defined as:

\begin{align}
	MLE = & \log(y^T\mathbf{K}^{-1}y) + \frac{1}{N}\log(det(\mathbf{K})) \nonumber \\
	= & \log(y^T\mathbf{K}^{-1}y) + \frac{1}{N}\sum_{i=1}^N{\log(\lambda_i(\mathbf{K}))}
	\label{mle_eq}
\end{align}

The right hand side of~\ref{mle_eq} can be broken into two terms: the first
term quantifies how well the model fits the data, which is simply a
Mahalanobis distance~\cite{mahalanobis1936generalized} between the model
predictions and the data; and the second term quantifies the model complexity
-- smoother covariance matrices will have smaller determinants and therefore
lower complexity penalties. Parameters of kernel function that minimize this
metric is the best set for the constructed surrogate. It should be noted here
that premature kernel parameter optimization can often lead to highly skewed
estimates and there should be performed only after sufficient number of data
points are evaluated~\cite{benassi2011robust, bull2011convergence}. A detailed discussion of hyper-parameter estimation in Bayesian optimization can be found
in~\cite{wang2014theoretical, jones1998efficient, shahriari2016taking}.

\subsection{Acquisition function}
\label{acqs}

\mydraft{
	\begin{itemize}
		\item selection of next point -- `infill' criterion
		\item list of commonly used acquisition functions , compare and contrast
		\item advantage of current implementation
	\end{itemize}
}

The acquisition function informs the selection of the next point for evaluation of the cost function, $y$. It is also sometimes called the 'in-fill' criterion.
The acquisition function is expected to locate the most useful point for evaluation of the cost function. This is particularly necessary when dealing with very
expensive cost functions. The acquisition function operates on the mean and
variance of the surrogate, $\tilde{y}$, and identifies the `best' point for next evaluation.
Typically, these functions are defined such that low acquisition values corresponds to
\begin{enumerate*}[label=\itshape\textbf{\alph*}\upshape)]
    \item low values of cost function (\emph{exploitation})
    \item high value of uncertainty (\emph{exploration}) 
    \item a balance between exploration and exploitation
\end{enumerate*}
Minimizing the acquisition function is used as a guide to select the next point of
evaluation. Some of the commonly used acquisition functions, that are implemented in the
software are given below:

\begin{enumerate}
	\item Lower Confidence Bound: Based on the mean and variance of the
		surrogate, the lower confidence bound of the surrogate is defined as
		\begin{equation}
			LCB(\mathbf{x}) = \mu(\mathbf{x}) - \kappa \sigma(\mathbf{x})
		\end{equation}
		$\kappa$ is a manually tunable parameter that signifies exploitative
		search for smaller values and explorative search for large values.
	\item Probability of Improvement: This estimates the probability of
		improvement of surrogate at a given location. It is given by:
		\begin{equation}
			PI(\mathbf{x}) = 
			\begin{cases}
				\phi(\frac{f_{min} - \mu(\mathbf{x})}{\sigma(\mathbf{x})}) , & \text{if } \sigma(\mathbf{x}) > 0 \\
				0, &\text{if } \sigma(\mathbf{x}) = 0
		\end{cases}
		\end{equation}
		where $$I(\mathbf{x}) = \text{min}(\frac{f_{min} -
		\mu(\mathbf{x})}{\sigma(\mathbf{x})}, 0)$$ is the effective improvement
		at a given location and $\phi(\mathbf{x})$ is the normal probability
		distribution function.
	\item Expected Improvement: This metric estimates the expectation of
		improvement of the surrogate at a location, given by:
		\begin{equation}
			EI(\mathbf{x}) = 
			\begin{cases}
				\mu(\mathbf{x}) + (f_{min} - \mu(\mathbf{x}))\Phi(\frac{f_{min} - 
				\mu(\mathbf{x})}{\sigma(\mathbf{x})}) + \sigma(\mathbf{x}) 
				\phi(\frac{f_{min}-\mu(\mathbf{x})}{\sigma(\mathbf{x})}),& 
				\text{if } \sigma(\mathbf{x}) > 0\\
				0,              & \text{if } \sigma(\mathbf{x}) = 0
			\end{cases}	
		\end{equation}
		Similar to lower confidence bound, expected improvement can also be
		parametrized for directing the optimization to exploration vs exploitation.
		This modified EI is called generalized EI (g-EI) and is very elaborately
		explained in~\cite{sasena2002exploration}.
\end{enumerate}

\subsection{Parallelization}
\label{parallel}
\mydraft{
	Multiple stages of parallelization possible:
	\begin{itemize}
		\item of linear algebra in surrogate construction
		\item selection of infill points -- multiple `surrogate' optimizers
		\item in terms of explicit and implicit parallelism
		\item function evaluations
			\begin{itemize}
				\item parallel execution of cost functions
				\item parallel evaluation of cost function -- multiple instances
			\end{itemize}
	\end{itemize}
}

There are several strategies for parallelization of this
algorithm~\cite{Wang2016a, Regis2007}. It could be parallelization of
\begin{enumerate*}[label=\itshape\textbf{\alph*}\upshape)]
	\item linear algebra during covariance matrix construction;
	\item acquisition function optimization for multiple optima/exploration-exploitation;
	\item cost function \emph{execution}; and
	\item cost function evaluation through multiple worker threads.
\end{enumerate*}
Of these, the first two strategies give little computational advantage as the
time spent in those stages is much smaller for expensive cost functions.
Therefore, we take the route of parallelizing cost function execution and
evaluation. This is especially critical and useful for HPC based evaluations,
which have the ability to run simultaneous multi-processor jobs. Also,
parallelizing function evaluation creates a natural setup for asynchronous data
assimilation. Since the update is Bayesian, the order in which the data is
added to the original prior is insignificant. The choice of next function
evaluation depends exclusively on the data and not the ordering of the data.
This can also be seen in the construction of the covariance/correlation matrix,
which is symmetric and the order of data only affects the order of rows and
columns and has no effect on the characteristics of the matrix. This property
of the algorithm enables asynchronocity to be built into it. The chosen
`in-fill' points for evaluation can be added later in the case of technical
delays (typically long queue times or hardware failures on a HPC cluster)
without affecting the progress of optimization. How much later and what are the
minimum number of evaluations per iteration form a part of analysis of this
algorithm which we present in the results section.

\section{Asynchronous Function Evaluation}
\label{async}

In order to abstract the cost function evaluation from the platform of
evaluation, we develop a cost function evaluator class. The \emph{evaluator} is
a special function that takes two parameters: a list of "new" points to
evaluate, and a list of "old" points to include if evaluation has completed.
The function evaluator returns three lists: a list of completed points (and
their cost function value), a list of pending points (that are still being
evaluated), and a list of points that failed to be evaluated. The union of the
two lists of input points is always equal to the union of the three lists of
output points.

The most straightforward implementation of a cost function evaluator is one
that loops through all input points and evaluates the cost function at each
point serially. We can extend this to evaluate points in parallel by spawning
worker threads (or processes) that work through the list of input points in
parallel until every point has been evaluated, then returning the results. In
both cases, the returned pending list is always empty: all points either
complete successfully or fail. We refer to both of these implementations as
"synchronous" because the Bayesian optimization algorithm does not continue to
the next iteration until all input points have been evaluated.

As cost function evaluation time may vary significantly, we introduce a
tuneable parameter that controls the fraction of ``new" points per iteration
that \emph{must complete} before the function evaluator returns. A value of
$1.0$ mimics synchronous/blocking behavior: all points must be evaluated before
the function evaluator returns, and the pending list will always be empty. A
value of 0.6 waits for at least 60\% of points from the current iteration to be
evaluated before continuing to the next iteration, returning the remaining
points in the pending list. A value of 0.0 will begin evaluating the input
points, then return immediately with all input points in the pending list. Only
the "new" points (the first argument of the function evaluator) are considered
part of this fraction to prevent a buildup of difficult-to-evaluate stale
points holding up an iteration.\footnote{This also prevents an iteration from
completing with mostly stale points.}

\begin{figure}
	\centering
	\includegraphics[width=0.5\textwidth]{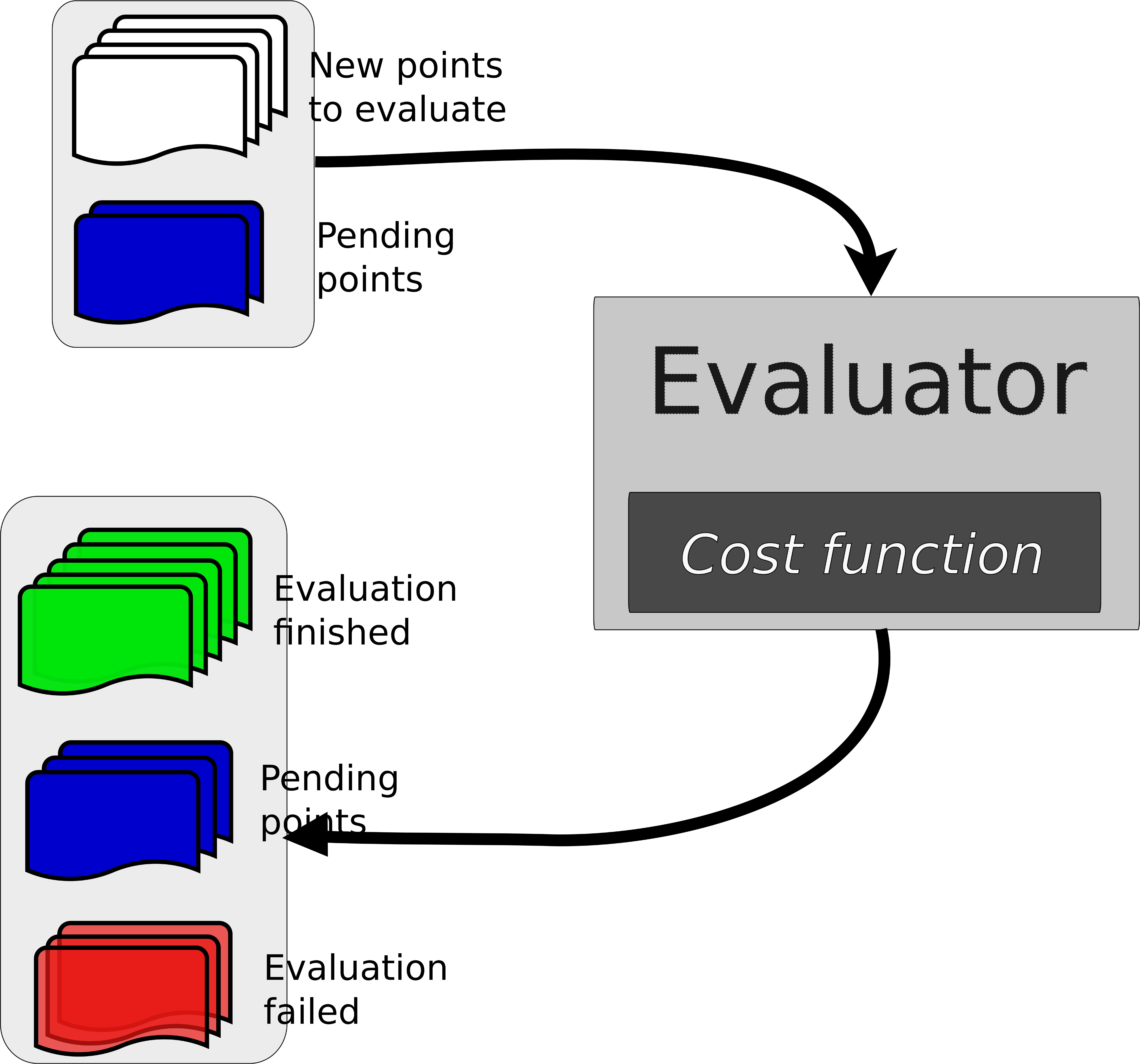}
	\caption{
		\textbf{Asynchronous evaluator}: The evaluator recieves two lists of
		pending points and new points for evaluation. It returns three lists --
		evaluation finished, pending points, and failed evaluation. The
		fraction of evaluations finished can be controlled by the fraction of
		asynchronocity (blocking-fraction).
	}
	\label{evaluator_class}
\end{figure}

Future calls to the function evaluator pass the previously-returned pending
points list as the second argument ("old" points). Once the evaluator finishes
evaluating the list of "new" points, it checks to see if any of the pending
"old" points have completed (or failed). If they have, they are added to the
appropriate list. If they have not, they are added to the pending list. The
evaluator never waits for "old" points.

\subsection{Implementation details} The package comes with an evaluator class
that supports asynchronous function evaluations. The basic inputs to an
asynchronous evaluator class are the \textit{maximum simultaneous} jobs
possible on the evaluation platform and the \textit{blocking-fraction} of jobs
to be finished before proceeding to the next bayesian update. For example, the
maximum jobs on a typical HPC cluster could be $50$ but on a local 4 core
machine is only $4$. The blocking fraction can vary from $0$ to $1$, with
\textit{zero} indicating fully asynchronous update while \textit{one} indicates
fully synchronous update. This class has to be derived for a specific platform,
examples are provided for a local machine and a SBATCH scheduler based HPC
system.
On a local machine, methods are implemented to spawn processes for each cost
function evaluation. Through the process id, the operating system's process
table is monitored to detect completion of each function evaluation. For a
cluster, sample evaluator sub-classes are included for an SBATCH based queing
system, for which methods are implemented to track the status of the job
through its \textit{jobId}. It comes with a complete setup for performing file
transfer and remote job submission, using Paramiko SSH library for Python. The
user only needs to provide external methods that can parse the output files
from a simulation and classify them into either of the three categories of:
\begin{enumerate*}[label=\textbf{\alph*})]
	\item \textit{ValueNotReady};
	\item \textit{EvaluationFailed}; and
	\item \textit{EvaluateAgain};
\end{enumerate*}
Jobs with \textit{ValueNotReady} are pushed into the pending list for future
Bayesian assimilation. If the evaluation is completed, the result parsing
script (user-provided) is invoked and it either returns a cost-function value
(float) or either of \textit{EvaluationFailed} or \textit{EvaluateAgain}. Jobs
with a \textit{EvaluationFailed} status will not be evaluated in the future.
Those with a \textit{EvaluateAgain} status will be submitted again.

\section{Results and discussion}
\label{results}
\mydraft{
	\begin{itemize}
		\item working of asynchronocity -- plot the number of pending jobs and
			completed jobs with iteration number
		\item show implementation enables multi-dimension optimization
	\end{itemize}
}

\subsection{Optimization}

\begin{figure}[h]
	\centering
	\begin{subfigure}{0.3\textwidth}
	\tiny
		\includegraphics[width=\textwidth, trim={0 0 0 35px}, clip]{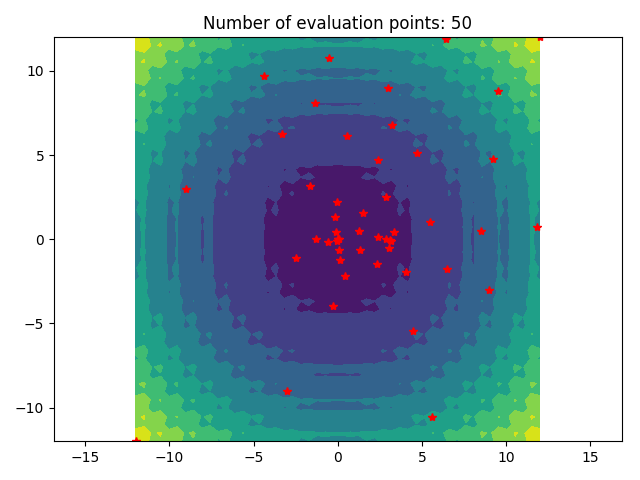}
		\caption{Function evaluations: Serial Bayesian optimization - only one
		evaluation per iteration}
	\end{subfigure}
	\begin{subfigure}{0.3\textwidth}
	\tiny
		\includegraphics[width=\textwidth, trim={0 0 0 35px}, clip]{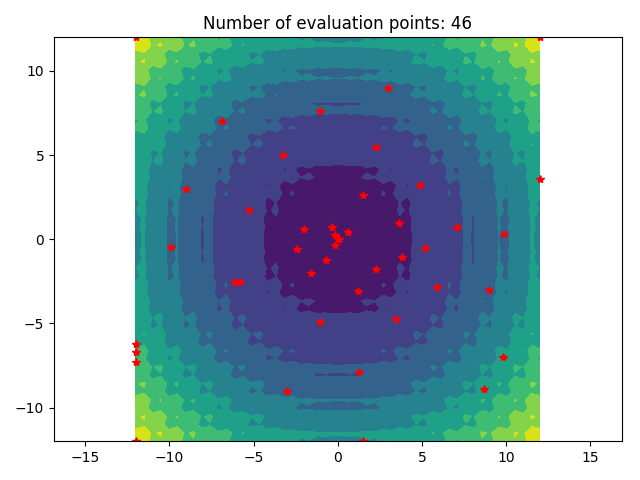}
		\caption{Function evaluations: Parallel Bayesian optimization with 4
		evaluations per iteration}
	\end{subfigure}
	\begin{subfigure}{0.3\textwidth}
	\tiny
		\includegraphics[width=\textwidth, trim={0 0 0 35px}, clip]{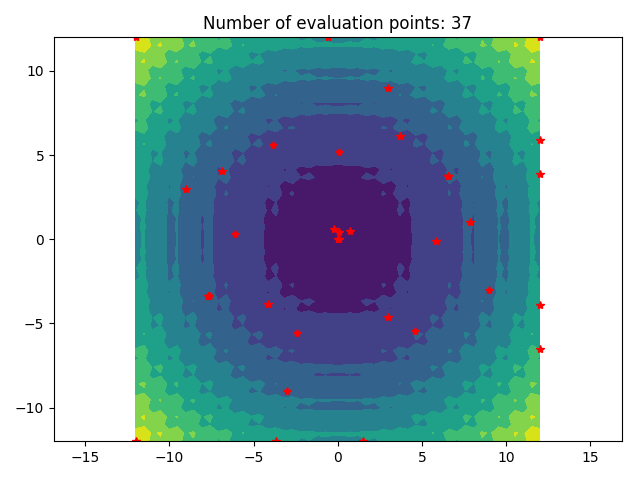}
		\caption{Function evaluations: Parallel Bayesian optimization with 8
		evaluations per iteration}
	\end{subfigure}
	\caption{Bayesian optimization on Rastrigin function with different levels
	of parallelization}
	\label{fig:serial_parallel}
\end{figure}

We first show that parallelization does not affect the ability of the algorithm
to find true global optima. Optimization was performed on the standard
Rastrigin [\ref{rastrigin}] cost function over a domain $(x,y) \in [-12,12]
\times [-12,12]$. The global minimum is located at $(0, 0)$ with a function
value of $0.0$. No asynchronocity was used and all the evaluations were used to
update the surrogate. In all the cases, a \textit{squared exponential} kernel
function was used with \textit{LCB} as the acquisition function. Across all
runs, the same $\kappa$ strategy was used for uniformity. Two different levels
of parallelization were used: the first with $4$ in-fill points per iteration
and the second with $8$ in-fill points per iteration. Both these cases were
compared with a standard serial Bayesian optimization. Figure
\ref{fig:serial_parallel} plots the locations of function evaluations with
different levels of parallelization. It should be observed that by using a
parallel optimizer, the optima is not compromised. The path of evaluation is
indeed affected by the degree of parallelization. The simultaneous
identification of in-fill points helps to mimic exploration-exploitation search
without explicitly tailoring the $\kappa$-strategy.

\begin{equation}
	f(x,y) = 20 + x^2 + y^2 - 10 cos(2\pi x) - 10 cos(2\pi y)
	\label{rastrigin}
\end{equation}

\subsubsection*{Asynchronous evaluator}

\begin{figure}[h]
	\centering
	\begin{subfigure}{0.45\textwidth}
		\raisebox{0px}{\includegraphics[width=\textwidth]{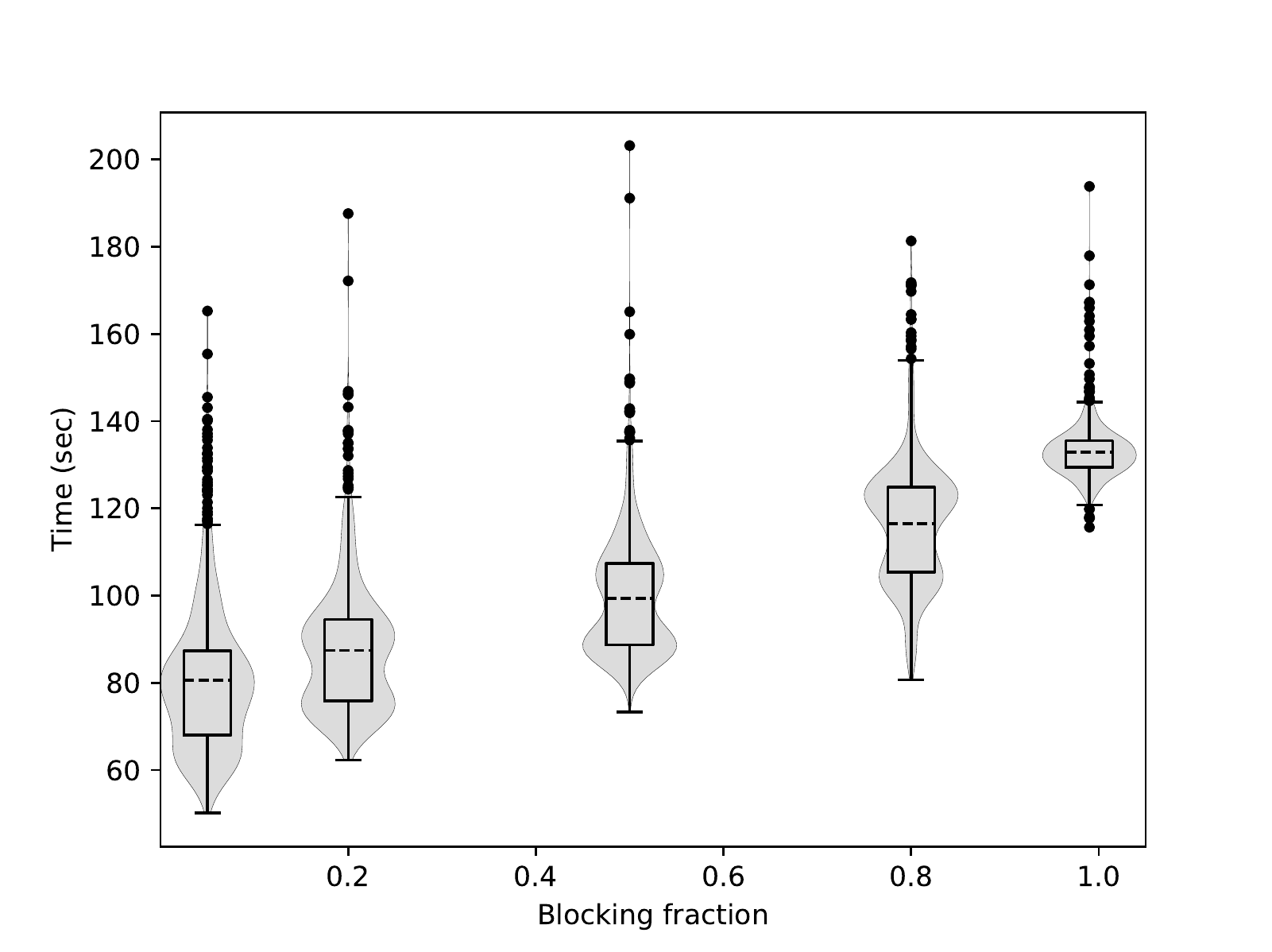}}
		\caption{Distribution of total optimization time for Bayesian optimization with asynchronous function evaluations}
		\label{fig:async_times_a}
	\end{subfigure}
	\begin{subfigure}{0.45\textwidth}
		\raisebox{13px}{\includegraphics[width=\textwidth]{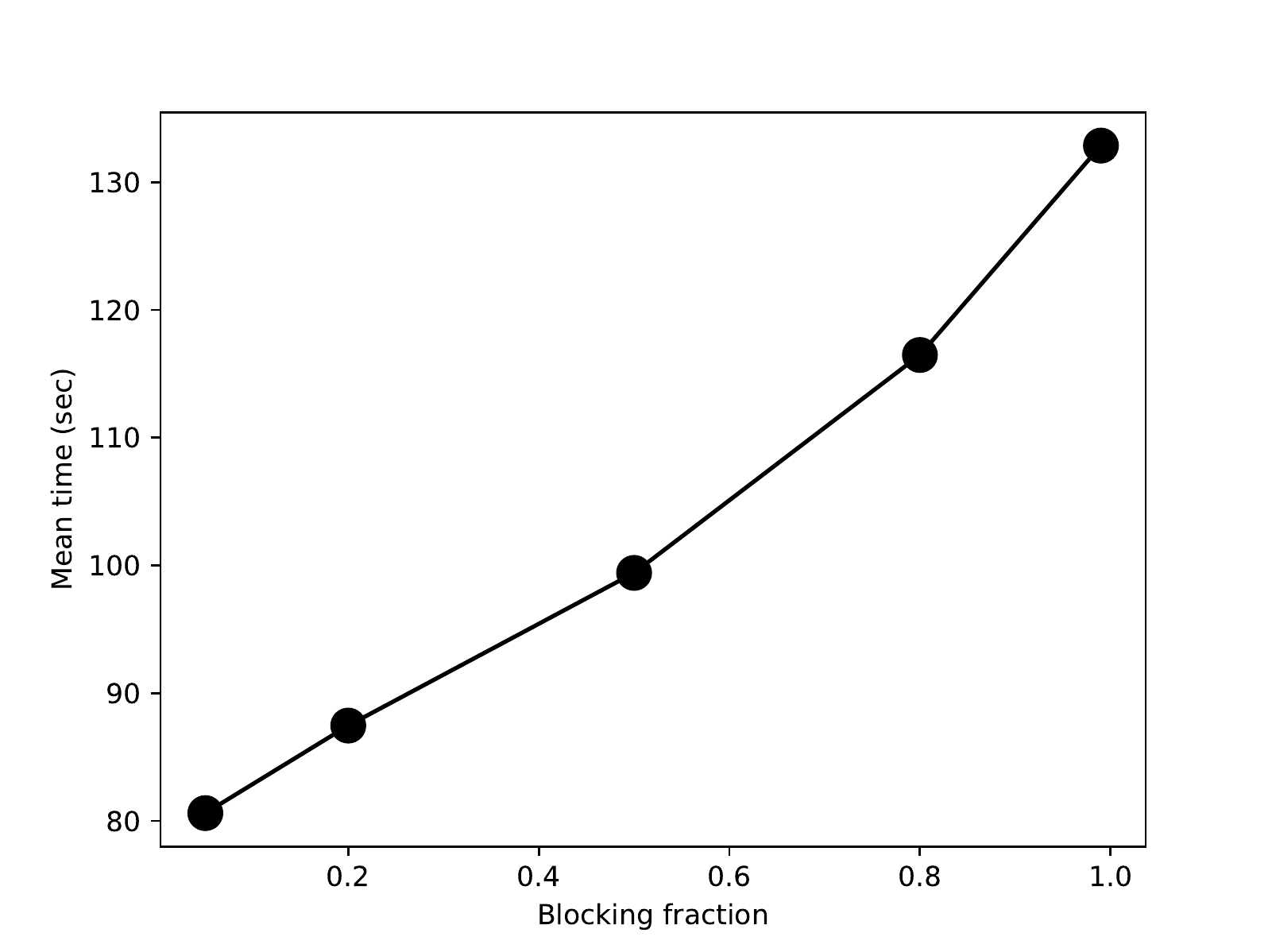}}
		\caption{Mean of completion times of Bayesian optimization with asynchronous function evaluations}
		\label{fig:async_times_b}
	\end{subfigure}
	\caption{Asynchronous Bayesian optimization : times of completion}
	\label{fig:async_times}
\end{figure}

Next, we demonstrate the performance and advantages of our asynchronous
evaluator. Using the same Rastrigin function [\ref{rastrigin}], we show how
asynchronous optimization can improve total optimization times under uncertain
function evaluation time. Since the cost function in this example is cheap,
function evaluations are instantaneous. Hence, to demonstrate asynchronous
working, we artificially enforce a finite total evaluation time. For this, we
chose to sample times from a normal distribution with a mean of $10$ sec and a
standard deviation of $2.5$ sec. This distribution is representative of the
queue wait time in standard schedulers. Other times, such as the time of
evaluation and time of startup are all assumed constant. The mean time of $10$
sec is chosen so that the total time is small and simultaneously not too small
to be overshadowed by program (Python, here) start up time. As with the
previous test, a \textit{squared exponential} kernel function was used along
with a \textit{LCB} acquisition function. Since the completion time (wait time,
here) of each evaluation is (deterministically) uncertain, several ($1000$)
realizations of the optimization were performed to get trends in total
optimization time. Figure~\ref{fig:async_times} shows the total time of
completion of the optimization for different values of blocking fractions.
Barring the few outliers (denoted by dots in figure~\ref{fig:async_times_a}),
it can be observed that there is a consistent reduction of total optimization
time with reducing blocking fraction. The mean completion time shows an
interesting linear trend with blocking fraction, indicating a strong advantage
of asynchronous optimization over a fully blocking optimization. In this simple
case, the improvement is upto $50\%$ on average. However, while there is a
significant reduction of the average time of completion, the uncertainty
increases for lower blocking fractions. This is an expected outcome of the
current test case, the asynchronous evaluator is highly influenced by the wait
times more than the fully blocking optimization. Finally, it should be noted
that in the worst case scenario, the total completion times for all blocking
fractions were similar and not affected by the presence of asynchronocity.

\subsubsection*{Benchmark Optimization}

Here, we show the working of the software with standard benchmark functions.
Table~\ref{tab:opt2},\ref{tab:opt5-10} shows the results. In all these cases, no parallelization
was used. The \textit{squared exponential} kernel was used with the
\textit{LCB} acquisition function. An annealing type $\kappa$
parameter~\cite{sasena2002exploration} was used for the LCB function. It can be
observed how PARyOpt can 

\begin{table}[H]
	\caption{Results from optimization of benchmark optimization test functions: 2D}
	\label{tab:opt2}
	\centering
	\begin{tabular}{|>{\centering\arraybackslash} m{5cm}|>{\centering\arraybackslash} m{3cm}|>{\centering\arraybackslash} m{4cm}|}
		\toprule
		Cost function & Optima location & Located Optima \\
		\midrule
		Ackley function & (0.0, 0.0) & (1.35e-11, -2.95e-08) \\
		\midrule
		Rosenbrock function &(1.0, 1.0) & (0.9999999945, 1.00000000124)\\
		\midrule
		Rastrigin function & (0.0, 0.0) & (1.97e-08, 1.17e-09)\\
		\midrule
		Griewangk function & (0.0, 0.0)& (-0.0003, -0.0009)\\
		\bottomrule
	\end{tabular}
\end{table}

\begin{table}[h]
	\caption{Results from optimization of benchmark optimization test functions: higher dimensions}
	\label{tab:opt5-10}
	\centering
	\begin{tabular}{|>{\centering\arraybackslash} m{3.5cm} >{\centering\arraybackslash} m{1cm}|>{\centering\arraybackslash} m{3.5cm}|>{\centering\arraybackslash} m{4cm}|}
		\toprule
		Cost function &  & Optima location & Distance of Located Optima from true optima ($l_\infty$)\\
		\midrule
		\multirow{3}{*}{Ackley function} &N=5D& \multirow{3}{*}{$x_i = 0.0$ for $i \in 1,2,..N$} & 1.75e-9 \\
										 &N=10D&  & 1.57e-7 \\
										 &N=20D&  & 1.57e-4 \\
		\midrule
		\multirow{3}{*}{Rosenbrock function} &N=5D&\multirow{3}{*}{$x_i = 1.0$ for $i \in 1,2,..N$}  & 1.25e-9 \\
											 &N=10D& & 3.62e-7 \\
											 &N=20D&  & 5.7e-4 \\
		\midrule
		\multirow{3}{*}{Rastrigin function} &N=5D&  \multirow{3}{*}{$x_i = 0.0$ for $i \in 1,2,..N$} & 1.136e-7\\
											&N=10D&  & 4.23e-6\\
											&N=20D&  & 3.8e-3 \\
		\midrule
		\multirow{3}{*}{Griewangk function} &N=5D& \multirow{3}{*}{$x_i = 0.0$ for $i \in 1,2,..N$} & 1.2e-4\\
											&N=10D& & 2.4e-3\\
											&N=20D&  & 5.9e-3 \\
		\bottomrule
	\end{tabular}
\end{table}

\subsection{Kriging}

Kriging is a more general interpolation technique based on Gaussian priors and
posteriors. It informs function values at various locations in terms of mean
and variance, and works similar to distance-weighted-interpolation. The
algorithmic basis of both Kriging and Bayesian optimization are very similar --
use local basis functions to create interpolant surrogate. In engineering,
there are several applications of kriging, especially in the realm of adaptive,
automated design of experiments. Due to the well known curse of dimensionality
and expensive experiments, engineers look to algorithms that inform regions of
maximum new information content. Our current implementation easily enables this
by selecting the exploratory mode of acquisition function.

We demonstrate this utility by considering the following example case: what is
the most efficient sampling for a given curve (1-D). We shall compare the
number of evaluations of Bayesian sampling to uniformly spaced evaluations,
typical of response surface reconstruction. We consider purely analytical
functions in this study -- although their evaluation is cheap, it provides an
easy powerful way to see the dependence of number of evaluations on the
complexity of the function.

Table \ref{tbl:krig1} shows how the current software can be used for kriging.
Each row shows varying modality of the underlying curve. The left column has
the result with unoptimized length parameters and the right column shows the
result with (MLE) optimized (kernel) parameters. Note the importance of
performing MLE minimization for achieving efficient reconstruction using
Kriging. All the results were done with a pure exploration mode using the \textit{LCB}
acquisition function.

\begin{table}[ht]
	\centering
	\begin{tabular}{|m{2.5cm}m{5.5cm}m{5.5cm}|}
		\toprule
		\centering
		Analytical function 	& Kriging without MLE optimization & Kriging with MLE optimization \\
		\midrule
		$y = x + 0.5 $		 & 
		\includegraphics[width=0.4\textwidth, trim={1.25cm .75cm 1.60cm 1.2cm}, clip]{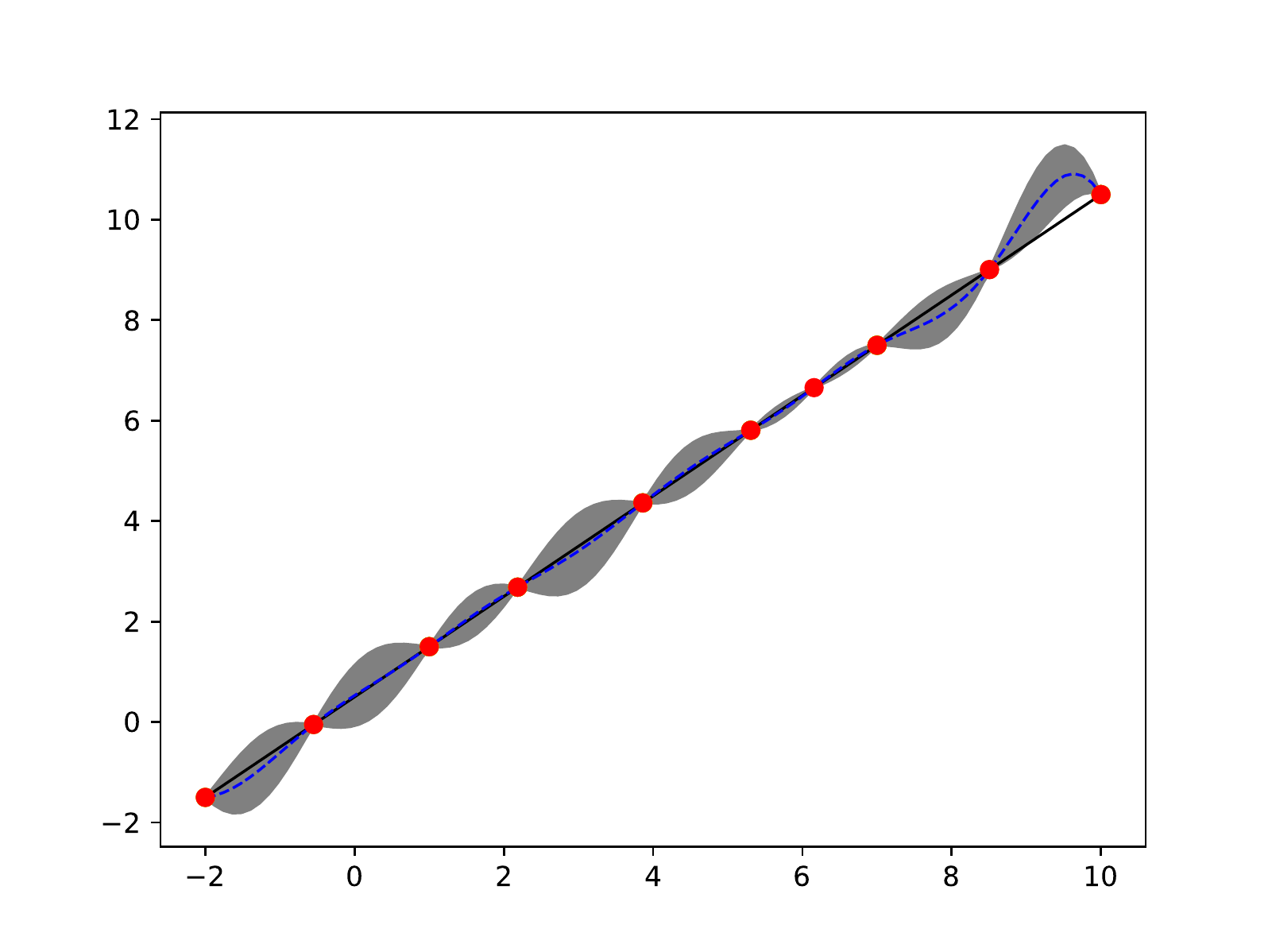}		&
		\includegraphics[width=0.4\textwidth, trim={1.25cm .75cm 1.60cm 1.2cm}, clip]{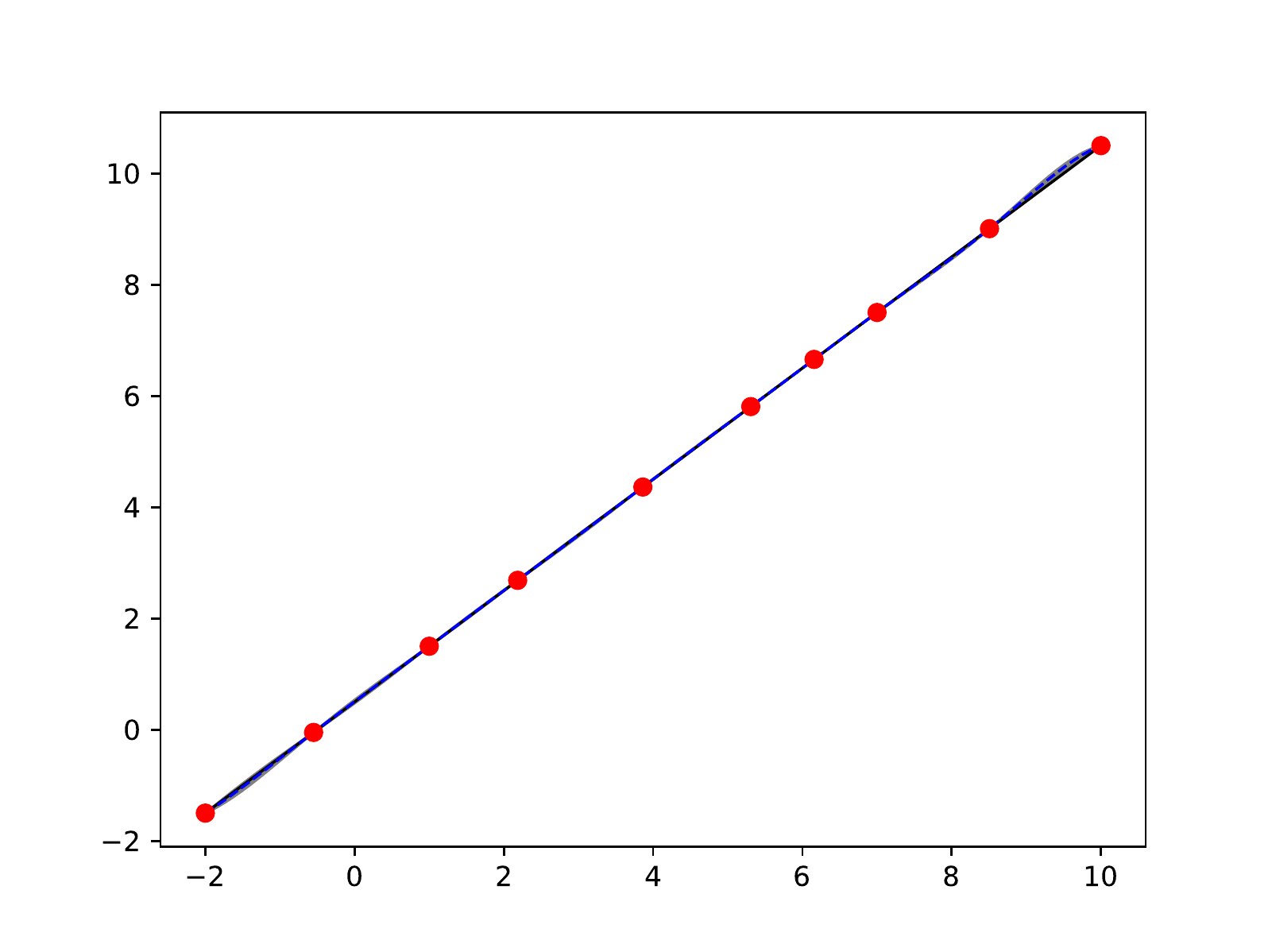} \\
		\midrule
		$y = (x-0.5)^2 + 1$ & 
		\includegraphics[width=0.4\textwidth, trim={1.25cm .75cm 1.60cm 1.2cm}, clip]{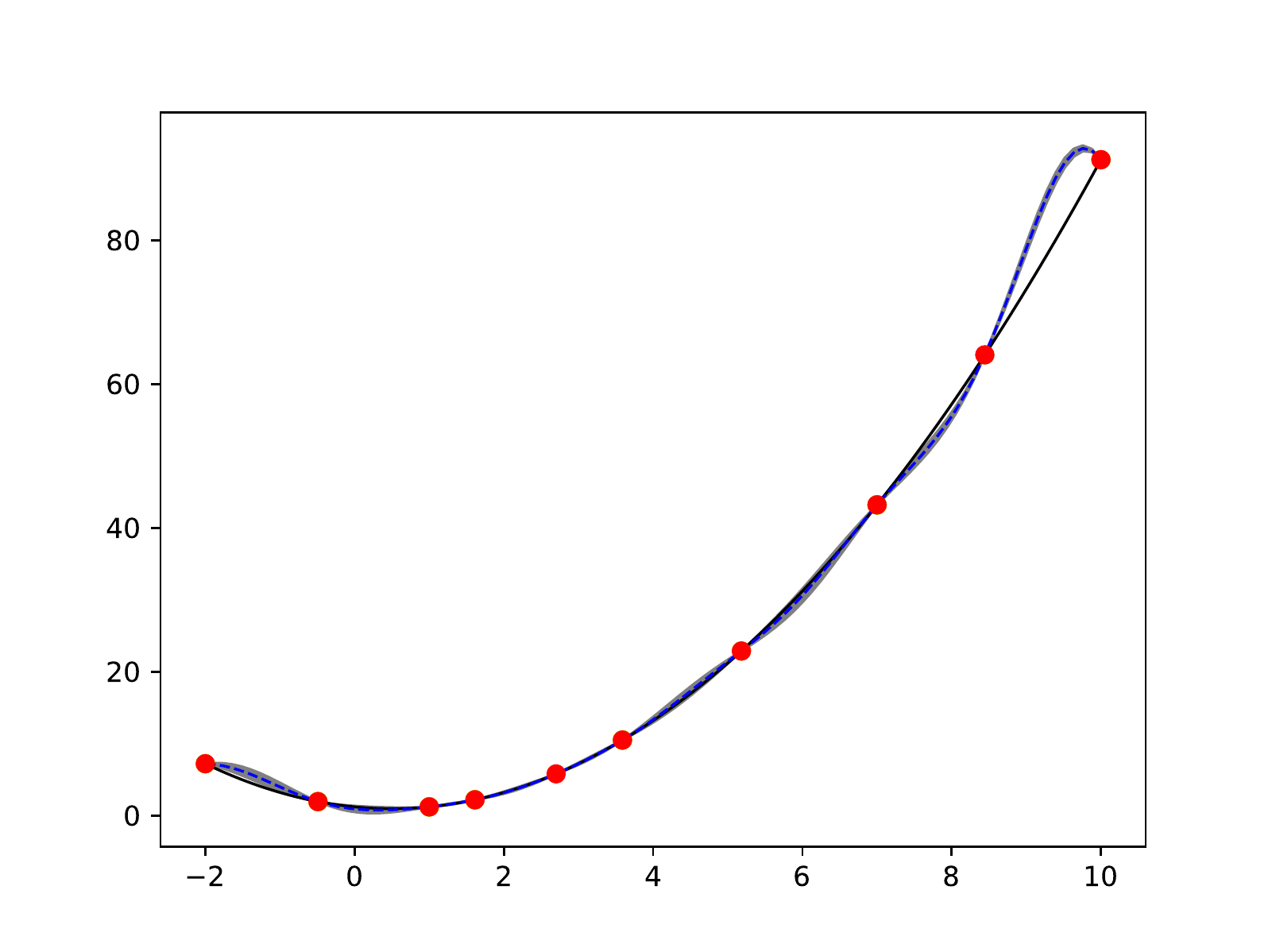}		&
		\includegraphics[width=0.4\textwidth, trim={1.25cm .75cm 1.60cm 1.2cm}, clip]{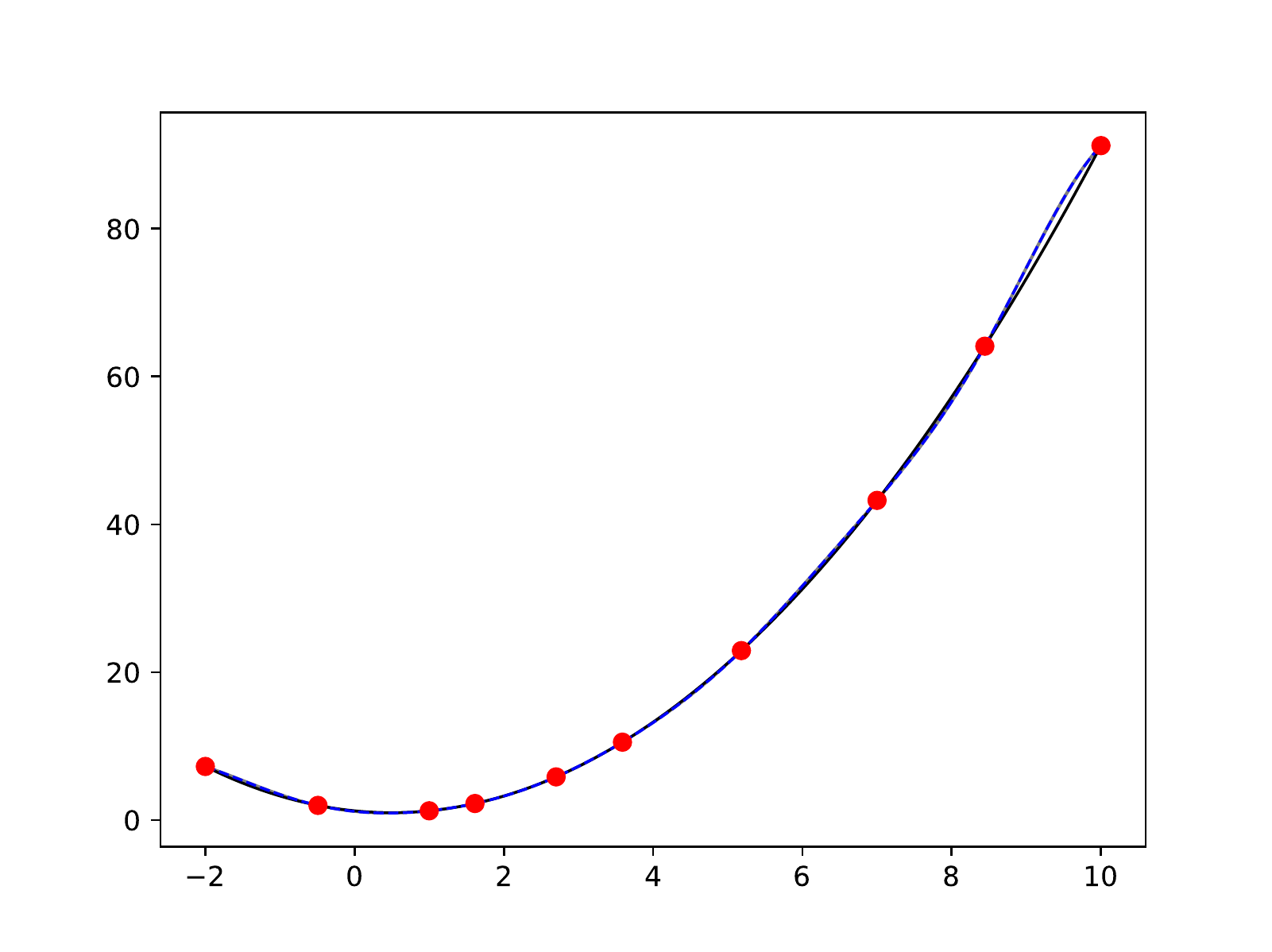} \\
		\midrule
		$y = sin(3 x^2 + (x-8)^2 + 1) $ & 
		\includegraphics[width=0.4\textwidth, trim={1.25cm .75cm 1.60cm 1.2cm}, clip]{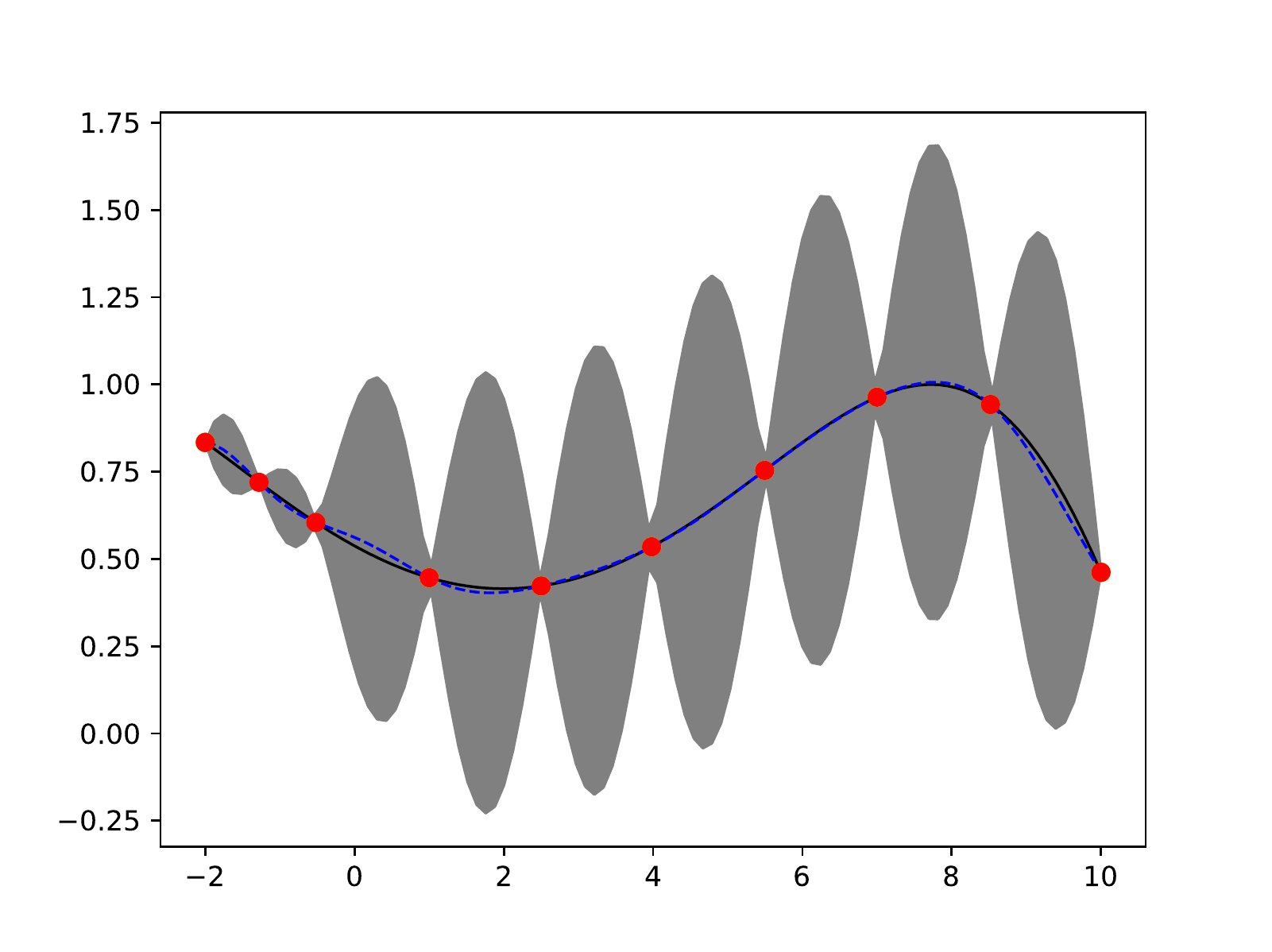}		&
		\includegraphics[width=0.4\textwidth, trim={1.25cm .75cm 1.60cm 1.2cm}, clip]{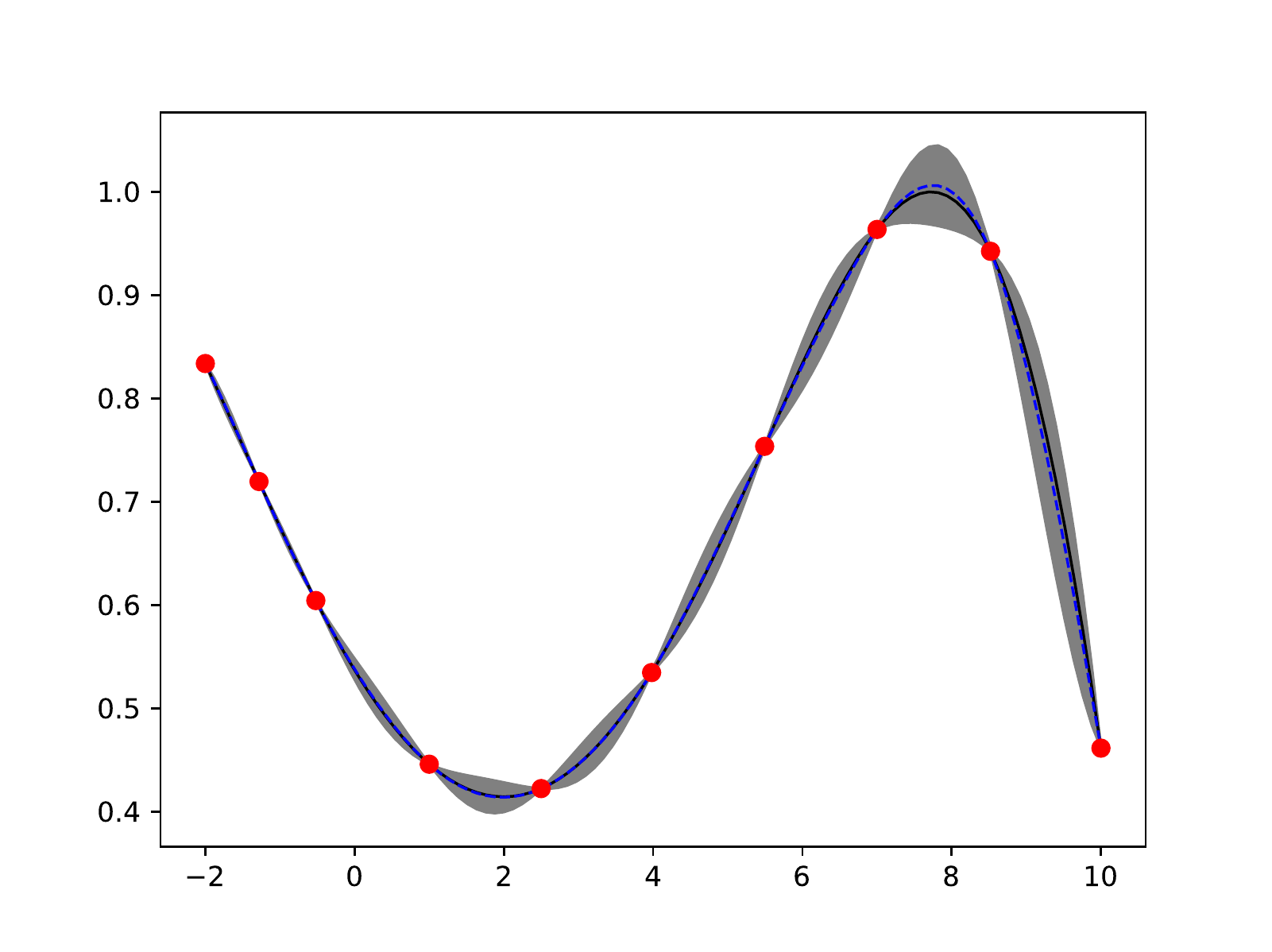} \\
		\midrule
		$y = sin(((x-6)/40)^2 + ((2x+1)/10)^3 ) $ & 
		\includegraphics[width=0.4\textwidth, trim={1.25cm .75cm 1.60cm 1.2cm}, clip]{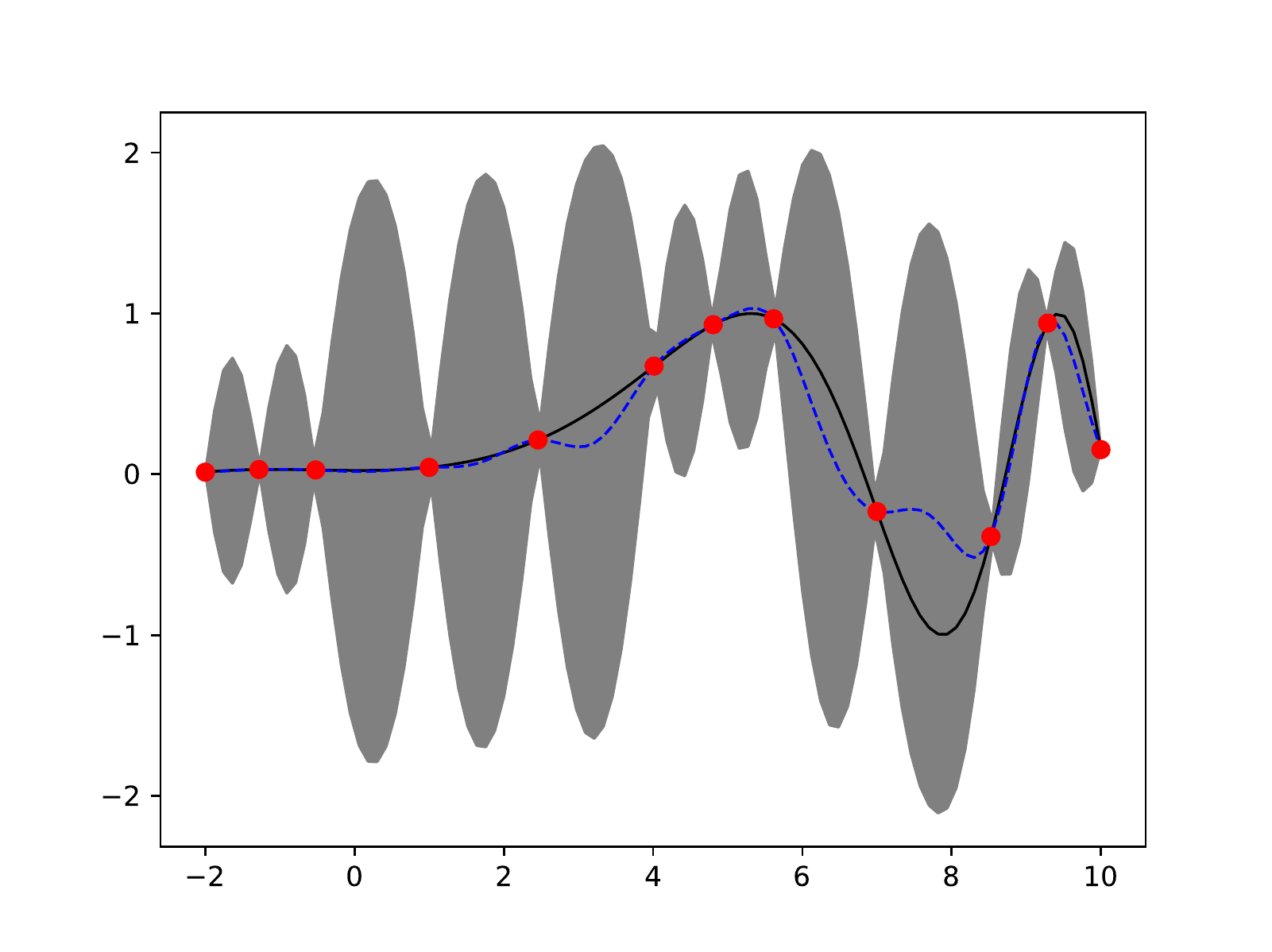}		& 
		\includegraphics[width=0.4\textwidth, trim={1.25cm .75cm 1.60cm 1.2cm}, clip]{./figs/3_kriging_example_optimized} \\
		\bottomrule
	\end{tabular}

	\caption{Kriging using PARyOpt. In these examples, we perform a budgeted
	Kriging, i.e., the total number of function evaluations are fixed. Note how
	hyper-parameter optimization using MLE helps to build greater confidence in the
	constructed surrogate.}

	\label{tbl:krig1}
\end{table}

\section{Conclusions and Future work}
\label{future}
In this work, we presented a software for performing parallelized asynchronous
Bayesian optimization. Several advantages of such asynchronous evaluations are
discussed and presented through test examples. The total time of optimization
improved up to $50\%$ using asynchronous parallel evaluations. This framework
is completely modular giving the user the ability to replace any functionality
with tailored functions. This framework is also extensible to perform Kriging. 
Examples were presented to show how the current framework can do this. This software is freely available on \href{www.bitbucket.org/baskargroup/paryopt}{Bitbucket} and
the documentation is hosted on \href{paryopt.rtfd.io}{ReadTheDocs}. We anticipate wide usage of this framework by the materials design community.

\mydraft{
	\begin{itemize}
		\item \sout{integration with other local optimizers - discuss about
				modularity in the code that enables integration with other
			optimization methods}
		\item \sout{hyper-parameter tuning}
		\item code availability
	\end{itemize}
}


\bibliographystyle{ACM-Reference-Format}

\bibliography{bo_bibliography}


\end{document}